\input amstex
\documentstyle{amsppt}

\magnification =\magstep 1

\topmatter
\title
On the twisted $q$-zeta functions and $q$-Bernoulli polynomials
\endtitle
\author {Taekyun Kim, Lee Chae Jang, Seog-Hoon Rim and Hong-Kyung Pak}
\endauthor
\affil{To the 62nd birthday of  Og-Yeon Choi}
\endaffil


\abstract{ One purpose of this paper is to define the twisted
$q$-Bernoulli numbers by using $p$-adic invariant integrals on
$\Bbb{Z}_p$. Finally, we construct the twisted $q$-zeta function
and $q$-$L$-series which interpolate the twisted $q$-Bernoulli
numbers.

}
\endabstract

\newline
\endkeywords
\endtopmatter

\document

\vskip 0.5 true cm {\bf 1. Introduction}
\par
\vskip 0.5 true cm

Throughout this paper $\Bbb{Z}_p, \Bbb{Q}_p, \Bbb{C}$ and
$\Bbb{C}_p$ are respectively denoted as the ring of $p$-adic
rational integers, the field of $p$-adic rational numbers, the
complex number field and the completion of algebraic closure of
$\Bbb{Q}_p$. The $p$-adic absolute value in $\Bbb{C}_p$ is
normalized so that $|p|_p = \frac 1p$. When one talks of
$q$-extension, $q$ is considered in many ways such as an
indeterminate, a complex number $q \in \Bbb{C}$, or a $p$-adic
number $q \in \Bbb{C}_p$. If $q \in \Bbb{C}$ one normally assumes
that $|q| < 1$. If $q \in \Bbb{C}_p$, we normally assume that
$|q-1|_p < p^{-1/(p - 1)}$ so that $q^x = \exp(x \log q)$ for
$|x|_p \le 1$. We use the notation as $$ [x] = [x : q] = 1 + q +
\cdots + q^{x-1},
$$ for $x \in \Bbb{Z}_p$.

Let $UD(\Bbb{Z}_p)$ be the set of uniformly differentiable
functions on $\Bbb{Z}_p$. For $f \in UD(\Bbb{Z}_p)$ the $p$-adic
$q$-integral was defined as
$$ I_q(f) = \int_{\Bbb{Z}_p} f(x) d \mu_q(x) = \lim_{N \rightarrow
\infty} \frac{1}{[p^N]} \sum_{0 \le x < p^N} f(x) q^x, \quad
\text{cf. [1, 2, 3]}. $$ Note that
$$ I_1(f) = \lim_{q \rightarrow 1} I_q(f) = \int_{\Bbb{Z}_p} f(x) d \mu_1(x)
= \lim_{N \rightarrow \infty} \frac{1}{p^N} \sum_{0 \le x < p^N}
f(x). $$

For a fixed positive integer $d$ with $(p,d) = 1$, let
$$  \aligned
X = X_d & = \lim_{\overleftarrow{N}} \Bbb{Z}/dp^N \Bbb{Z}, \quad X_1 = \Bbb{Z}_p, \\
X^* & = \cup_{\ssize{0 < a < dp} \atop \ssize{(a,p)=1}} (a + dp
\Bbb{Z}_p), \\
a & + dp^N \Bbb{Z}_p = \{ x \in X | x \equiv a (\text{mod} \ dp^N)
\},
\endaligned $$
where $a \in \Bbb{Z}$ lies in $0 \le a < dp^N$, cf. [1, 2, 3, 4,
5, 6]. Let $\Bbb{Z}$ be the set of integers. For $h  \in \Bbb{Z},$
$k\Bbb \in \Bbb N,$ the $q$-Bernoulli polynomials were defined as
$$ \beta_n^{(h,k)}(x,q) =  \int_{\Bbb{Z}_p^k} [x+x_1 + \cdots + x_k
]^n q^{\overset k \to {\underset {i=1} \to \sum} x_i (h-i)}
d\mu_q(x_1) \cdots d\mu_q(x_k), \quad \text{cf. [1, 2]}. \tag 1 $$
The $q$-Bernoulli polynomials at $x=0$ are called $q$-Bernoulli
numbers. In [1] it was shown that the $q$-Bernoulli numbers were
written as
$$ \beta_n^{(h,k)} (= \beta_n^{(h,k)}(q)) = \beta_n^{(h,k)}(0,q).
$$
Indeed, $\lim_{q\rightarrow 1} \beta_n^{(h,k)}(q) = B_n^{(k)}$,
where $B_n^{(k)}$ are Bernoulli numbers of order $k,$ see [1, 2,
3]. Let $\chi$ be a Dirichlet character with conductor $f\in\Bbb N
$. Then the Dirichlet $L$-series attached to $\chi$ is defined as
$$ L(s,\chi) = \sum_{n=1}^\infty \frac{\chi(n)}{n^s}, \text{ for $ s \in
\Bbb{C},$ cf. [7, 8]}.$$ When $\chi=1$, this is the Riemann zeta
function. In [1], $q$-analogue of $\zeta$-function was defined as
follows: For $h \in \Bbb{Z}, s \in \Bbb{C},$
$$ \zeta_q^{(h)}(s,x) = \frac{1-s+h}{1-s} (q - 1) \sum_{n=0}^\infty
\frac{q^{(n+x)h}}{[n+x]^{s-1}} +  \sum_{n=0}^\infty
\frac{q^{(n+x)h}}{[n+x]^{s}}. \tag 2 $$ Note that
$\zeta_q^{(h)}(s,x)$ is an analytic continuation on $\Bbb{C}$
except for $s=1$ with
$$ \zeta_q^{(h)}(1-m,x) = -\frac{\beta_m^{(h,1)}(x,q)}{m}, \text{ for $ m
\in\Bbb N.$} \tag 3 $$ In [1], we easily see that
$$ \beta_m^{(h,1)}(x,q) = -m \sum_{n=0}^\infty
q^{(n+x)h} [n+x]^{m-1} - (q-1)(m+h)  \sum_{n=0}^\infty q^{(n+x)h}
[n+x]^{m}. \tag 4 $$ It follows from (2) that
$$ \lim_{q \rightarrow 1} \zeta_q^{(h)}(s,x) = \zeta(s,x) =
\sum_{n=0}^\infty \frac{1}{(n+x)^s}. $$ By the meaning of the
$q$-analogue of Dirichlet $L$-series, we consider the following
$L$-series:
$$ L_q^{(h)}(s,\chi) = \frac{1-s+h}{1-s} (q - 1) \sum_{n=1}^\infty
\frac{q^{nh} \chi(n)}{[n]^{s-1}} +  \sum_{n=1}^\infty \frac{q^{nh}
\chi(n)}{[n]^{s}}, \text{ cf. [1], } \tag 5 $$ for $h \in \Bbb{Z},
s \in \Bbb{C}$. It is easy to see that $L_q^{(h)}(s,\chi)$ is an
analytic continuation on $\Bbb{C}$ except for $s=1$. For $m \geq
0,$ the generalized extended $q$-Bernoulli numbers with $\chi$ are
defined as
$$ \beta_{m,\chi}^{(h,k)}(q) = [f]^{m-k} \sum_{i_1, \cdots, i_k = 0}^{f-1}
q^{\overset k \to {\underset {l=1} \to \sum} (h-l+1)i_l}
\beta_{m}^{(h,k)} (\frac{\overset k \to {\underset {l=1} \to \sum}
i_l}{f}, q^f)\left(\overset k \to {\underset {j=1} \to \Pi}
\chi(i_j)\right), \text{ see [1] }. \tag 6 $$  By (5) and (6), we
easily see that
$$ L_q^{(h)}(1-m,\chi)= - \frac{\beta_{m,\chi}^{(h,1)}}{m},\text{ for $m\in \Bbb N$, cf. [1]}. $$

In the present paper we give  twisted $q$-Bernoulli numbers by
using $p$-adic invariant integrals on $\Bbb{Z}_p$. Moreover, we
construct the analogs of $q$-zeta function and $q$-$L$-series
which interpolate the twisted $q$-Bernoulli numbers at negative
integers.


\vskip 0.5 true cm {\bf 2. $q$-extension of Bernoulli numbers}
\par
\vskip 0.5 true cm

In this section we assume that $q \in \Bbb{C}_p$ with $|1- q|_p <
1$. By the definition  of $p$-adic invariant integrals, we see
that
$$ I_1(f_1) = I_1(f) + f^{\prime}(x), \tag 7$$
where $f_1(x) = f(x+1).$ Let
$$T_p = \cup_{n \ge 1}C_{p^n} = \lim_{n\rightarrow \infty}
C_{p^n}, $$ where $C_{p^n} = \{ w | w^{p^n} = 1 \}$ is the cyclic
group of order $p^n$. For $w \in T_p ,$ we denote by $\phi_w :
\Bbb{Z}_p \longrightarrow \Bbb{C}_p$ the locally constant function
$x \mapsto w^x$. If we take $f(x) = \phi_w(x) e^{tx}$, then we
easily see that
$$ \int_{\Bbb{Z}_p} e^{tx} \phi_w(x) d\mu_1(x) = \frac{t}{we^t - 1}, \quad
\text{cf. [5]}. \tag 8 $$ It is obvious from (7) that
$$ \int_{X} e^{tx} \chi(x) \phi_w(x) d\mu_1(x) =
\frac{\overset f \to {\underset {i=1} \to \sum} \chi(i) \phi_w(i)
e^{it}} {w^f e^{ft} - 1}, \text{ cf. [5] }. \tag 9 $$ Now we
define the analogue of Bernoulli numbers as follows:
$$ \aligned
e^{xt}  \frac{t}{we^t - 1} & = \sum_{n=0}^\infty B_{n,w}(x) \frac{t^n}{n!} ,\\
\frac{\overset f \to {\underset {i=1} \to \sum} \chi(i) \phi_w(i)
e^{it}} {w^f e^{ft} - 1} & = \sum_{n=0}^\infty B_{n,w,\chi}
\frac{t^n}{n!}, \text{ cf. [5] }.
\endaligned \tag 10 $$
By (8), (9) and (10), it is not difficult to see that
$$ \int_{\Bbb{Z}_p} x^n \phi_w(x) d\mu_1(x) = B_{n,w} \quad \text{and} \quad
 \int_{X} \chi(x) x^n \phi_w(x) d\mu_1(x) = B_{n,w,\chi} , \tag 11$$
 where $B_{n,w} = B_{n,w}(0)$.

 From (11) we consider twisted $q$-Bernoulli numbers
  using $p$-adic $q$-integral on $\Bbb{Z}_p$. For $w \in T_p$ and $h \in \Bbb{Z},$
 we define the twisted $q$-Bernoulli polynomials as
 $$ \beta^{(h)}_{m,w}(x,q) = \int_{\Bbb{Z}_p} q^{(h-1)y} w^y [x+y]^m d\mu_q(y).
 \tag 12 $$
 Observe that
 $$ \lim_{q\rightarrow 1}  \beta^{(h)}_{m,w}(x, q) = B_{m,w}(x). $$
When $x=0$, we write  $\beta^{(h)}_{m,w}(0,q) =
\beta^{(h)}_{m,w}(q)$, which are called twisted $q$-Bernoulli
numbers. It follows from (12) that
$$ \beta^{(h)}_{m,w}(x,q) = \frac{1}{(1-q)^{m-1}} \sum_{k=0}^m {m \choose k} q^{xk}
(-1)^k \frac{k+h}{1 - q^{h+k} w}. \tag 13 $$ The Eq.(13) is
equivalent to
$$\beta^{(h)}_{m,w}(q) = -m \sum_{n=0}^\infty [n]^{m-1}q^{hn} w^n -
(q-1)(m+h) \sum_{n=0}^\infty [n]^{m}q^{hn} w^n. \tag 14 $$ From
(13), we obtain the below distribution relation for the twisted
$q$-Bernoulli polynomials as follows: For $n \ge 0,$
$$ \beta^{(h)}_{n,w}(x,q) = [f]^{n-1} \sum_{a=0}^{f-1} w^a q^{ha} \beta^{(h)}_{n,w^f}(\frac
{a+x}f,q^f). $$ Let $\chi$ be the Dirichlet character with
conductor $f \in \Bbb N$. Then we define the generalized twisted
$q$-Bernoulli numbers as follows: For $n \ge 0$,
$$ \beta^{(h)}_{m,w,\chi}(q) = \int_X \chi(x) q^{(h-1)x}w^x [x]^m
d\mu_q(x). \tag 15 $$  By (15), it is easy to see that
$$ \beta^{(h)}_{m,w,\chi}(q) = [f]^{m-1} \sum_{a=0}^{f-1} \chi(a) w^a q^{ha}
\beta^{(h)}_{m,w^f}(\frac af,q^f). \tag 16 $$  Remark. We note
that $\lim_{q\rightarrow 1} \beta^{(h)}_{m,w,\chi}(q) =
B_{m,w,\chi},$ ( see Eq. (10) ).


\vskip 0.5 true cm {\bf 3. $q$-zeta functions}
\par
\vskip 0.5 true cm

In this section we assume that $q \in \Bbb{C}$ with $|q| < 1$.
Here we construct the twisted $q$-zeta function
 and the twisted $q$-$L$-series (see Eq.(2) and Eq.(3)). Let $\Bbb{R}$ be the field of
real numbers and let $w$ be the $p^r$-th root of unity. For $q \in
\Bbb{R}$ with $0 < q < 1, s \in \Bbb{C}$ and $h \in \Bbb{Z},$ we
define the twisted $q$-zeta function as
$$ \zeta_{q,w}^{(h)}(s) = \frac{1 - s + h}{1 - s} (q-1)
\sum_{n=1}^\infty \frac{q^{nh} w^n}{[n]^{s-1}} + \sum_{n=1}^\infty
\frac{q^{nh} w^n}{[n]^{s}}. \tag 17 $$ Note that
$\zeta_{q,w}^{(h)}(s)$ is an analytic continuation on $\Bbb{C}$
except for $s=1$ and $\underset {q\rightarrow 1} \to \lim
\zeta^{(h)}_{q,w}(s) =
\zeta(s,w)=\sum_{n=1}^{\infty}\frac{w^n}{n^s},$ cf. [4].
 We see,
by (17), that
$$ \zeta^{(s-1)}_{q,w}(s) = \sum_{n=1}^\infty \frac{q^{n(s-1)}
w^n}{[n]^{s}}. \tag 18 $$

In what follows, the notation $ \zeta^{(s-1)}_{q,w}(s)$ will be
replaced by  $\zeta_{q,w}(s)$, that is,
$$ \zeta_{q,w}(s) (=  \zeta^{(s-1)}_{q,w}(s))= \sum_{n=1}^\infty \frac{q^{n(s-1)}
w^n}{[n]^{s}}. $$ We note that Eq.(18) is the $q$-extension of
Riemann zeta function. By (14) and (17) we give the values of
$\zeta^{(h)}_{q,w}(s)$ at negative integers as follows: For $m \in
\Bbb N,$ we have
$$ \zeta^{(h)}_{q,w}(1-m) = - \frac{\beta^{(h)}_{m,w}(q)}{m}. \tag
19 $$ By (17), we also see that
$$ \zeta_{q,w}(1-m) = \sum_{n=1}^\infty [n]^{m-1} q^{-mn} w^n . \tag 20 $$
The Eq.(20) seems to be the $q$-analogue of  Euler divergence
theorem for Riemann zeta function. Now we also consider the
twisted $q$-analogue of Hurwitz zeta function as follows: For $s
\in \Bbb{C},$ define
$$ \zeta^{(h)}_{q,w}(s,x) = \frac{1 - s + h}{1 - s} (q-1)
\sum_{n=0}^\infty \frac{q^{(n+x)h} w^n}{[n+x]^{s-1}} +
\sum_{n=0}^\infty \frac{q^{(n+x)h} w^n}{[n+x]^{s}}. \tag 21 $$

 Note that $\zeta_{q}^{(h)}(s,x)$ has an analytic
continuation on $\Bbb{C}$ with only one simple poles at  $s=1$. By
Eq.(13), Eq.(14) and Eq.(21), we obtain the following:
$$  \zeta^{(h)}_{q,w}(1-m,x) = - \frac{\beta^{(h)}_{m,w}(x,q)}{m},
\quad \text{for} \ m > 0. $$

Now we consider the twisted $q$-$L$-series which interpolate
twisted generalized $q$-Bernoulli numbers as follows: For $s \in
\Bbb{C},$ define
$$ L^{(h)}_{q,w}(s, \chi) = \frac{1 - s + h}{1 - s} (q-1)
\sum_{n=1}^\infty \frac{q^{nh} w^n \chi(n)}{[n]^{s-1}} +
\sum_{n=1}^\infty \frac{q^{nh} w^n \chi(n)}{[n]^{s}}, \tag 22 $$
where $w$ is the $p^r$-th root of unity.

For any positive integer $m$ we have
$$ L^{(h)}_{q,w}(1-m, \chi) = - \frac{\beta^{(h)}_{m,w,\chi}(q)}{m}.
$$
The Eq.(22) implies that
$$ \aligned
L^{(s-1)}_{q,w}(s, \chi) & = \sum_{n=1}^\infty \frac{q^{n(s-1)}
w^n
\chi(n)}{[n]^{s}} \\
& = [f]^{-s} \sum_{a=1}^{f} \chi(a) w^a q^{(s-1)a}
\zeta_{q^f,w^f}(s,\frac af). \endaligned $$

\noindent {\bf Question.} Find a $q$-analogue of the $p$-adic
twisted $L$-function which interpolates $q$-Bernoulli numbers
$\beta^{(h)}_{m,w,\chi}(q)$, cf. [3].

ACKNOWLEDGEMENTS: This paper was supported by Korea Research
Foundation Grant( KRF-2002-050-C00001).


\vskip0.5cm

\noindent Institute of Science Education, Kongju National
University, Kongju 314-701, Korea, tkim\@kongju.ac.kr \par
\noindent Department of Mathematics and Computer Science, KonKuk
University, Choongju, Chungbuk 380-701, Korea,
leechae.jang\@kku.ac.kr \par \noindent Department of Mathematics
Education, Kyungpook University, Taegu,
 702-701,  Korea, shrim\@kyungpook.ac.kr\par
 \noindent Faculty of Informatoin and
Science, Daegu Haany University, Kyungsan,  712-240, Korea,
hkpak\@ik.ac.kr

\enddocument

-----------------------------7d530c114a4
Content-Disposition: form-data; name="archive"

math
-----------------------------7d530c114a4
Content-Disposition: form-data; name="is_an_author"

1
-----------------------------7d530c114a4
Content-Disposition: form-data; name="Title"

On the twisted $q$-zeta functions and $q$-Bernoulli polynomials
-----------------------------7d530c114a4
Content-Disposition: form-data; name="Authors"

Taekyun Kim, L.C. Jang, S.H.Rim, H.K. Pak
-----------------------------7d530c114a4
Content-Disposition: form-data; name="Subj-class"

NT
-----------------------------7d530c114a4
Content-Disposition: form-data; name="MSC-class"

11B68, 11M06
-----------------------------7d530c114a4
Content-Disposition: form-data; name="Comments"

6 pages
-----------------------------7d530c114a4
Content-Disposition: form-data; name="Report-no"

-----------------------------7d530c114a4
Content-Disposition: form-data; name="Journal-ref"

-----------------------------7d530c114a4
Content-Disposition: form-data; name="DOI"

-----------------------------7d530c114a4
Content-Disposition: form-data; name="Abstract"

We study the twisted q-zeta functions and twisted q-Bernoulli polynomials
-----------------------------7d530c114a4
Content-Disposition: form-data; name="f1"; filename="C:\WINDOWS\¹ÙÅÁ È­¸é\twisted.tex"
Content-Type: application/octet-stream

\input amstex
\documentstyle{amsppt}

\magnification =\magstep 1

\topmatter
\title
On the twisted $q$-zeta functions and $q$-Bernoulli polynomials
\endtitle
\author {Taekyun Kim, Lee Chae Jang, Seog-Hoon Rim and Hong-Kyung Pak}
\endauthor
\affil{To the 62nd birthday of  Og-Yeon Choi}
\endaffil


\abstract{ One purpose of this paper is to define the twisted
$q$-Bernoulli numbers by using $p$-adic invariant integrals on
$\Bbb{Z}_p$. Finally, we construct the twisted $q$-zeta function
and $q$-$L$-series which interpolate the twisted $q$-Bernoulli
numbers.

}
\endabstract

\newline
\endkeywords
\endtopmatter

\document

\vskip 0.5 true cm {\bf 1. Introduction}
\par
\vskip 0.5 true cm

Throughout this paper $\Bbb{Z}_p, \Bbb{Q}_p, \Bbb{C}$ and
$\Bbb{C}_p$ are respectively denoted as the ring of $p$-adic
rational integers, the field of $p$-adic rational numbers, the
complex number field and the completion of algebraic closure of
$\Bbb{Q}_p$. The $p$-adic absolute value in $\Bbb{C}_p$ is
normalized so that $|p|_p = \frac 1p$. When one talks of
$q$-extension, $q$ is considered in many ways such as an
indeterminate, a complex number $q \in \Bbb{C}$, or a $p$-adic
number $q \in \Bbb{C}_p$. If $q \in \Bbb{C}$ one normally assumes
that $|q| < 1$. If $q \in \Bbb{C}_p$, we normally assume that
$|q-1|_p < p^{-1/(p - 1)}$ so that $q^x = \exp(x \log q)$ for
$|x|_p \le 1$. We use the notation as $$ [x] = [x : q] = 1 + q +
\cdots + q^{x-1},
$$ for $x \in \Bbb{Z}_p$.

Let $UD(\Bbb{Z}_p)$ be the set of uniformly differentiable
functions on $\Bbb{Z}_p$. For $f \in UD(\Bbb{Z}_p)$ the $p$-adic
$q$-integral was defined as
$$ I_q(f) = \int_{\Bbb{Z}_p} f(x) d \mu_q(x) = \lim_{N \rightarrow
\infty} \frac{1}{[p^N]} \sum_{0 \le x < p^N} f(x) q^x, \quad
\text{cf. [1, 2, 3]}. $$ Note that
$$ I_1(f) = \lim_{q \rightarrow 1} I_q(f) = \int_{\Bbb{Z}_p} f(x) d \mu_1(x)
= \lim_{N \rightarrow \infty} \frac{1}{p^N} \sum_{0 \le x < p^N}
f(x). $$

For a fixed positive integer $d$ with $(p,d) = 1$, let
$$  \aligned
X = X_d & = \lim_{\overleftarrow{N}} \Bbb{Z}/dp^N \Bbb{Z}, \quad X_1 = \Bbb{Z}_p, \\
X^* & = \cup_{\ssize{0 < a < dp} \atop \ssize{(a,p)=1}} (a + dp
\Bbb{Z}_p), \\
a & + dp^N \Bbb{Z}_p = \{ x \in X | x \equiv a (\text{mod} \ dp^N)
\},
\endaligned $$
where $a \in \Bbb{Z}$ lies in $0 \le a < dp^N$, cf. [1, 2, 3, 4,
5, 6]. Let $\Bbb{Z}$ be the set of integers. For $h  \in \Bbb{Z},$
$k\Bbb \in \Bbb N,$ the $q$-Bernoulli polynomials were defined as
$$ \beta_n^{(h,k)}(x,q) =  \int_{\Bbb{Z}_p^k} [x+x_1 + \cdots + x_k
]^n q^{\overset k \to {\underset {i=1} \to \sum} x_i (h-i)}
d\mu_q(x_1) \cdots d\mu_q(x_k), \quad \text{cf. [1, 2]}. \tag 1 $$
The $q$-Bernoulli polynomials at $x=0$ are called $q$-Bernoulli
numbers. In [1] it was shown that the $q$-Bernoulli numbers were
written as
$$ \beta_n^{(h,k)} (= \beta_n^{(h,k)}(q)) = \beta_n^{(h,k)}(0,q).
$$
Indeed, $\lim_{q\rightarrow 1} \beta_n^{(h,k)}(q) = B_n^{(k)}$,
where $B_n^{(k)}$ are Bernoulli numbers of order $k,$ see [1, 2,
3]. Let $\chi$ be a Dirichlet character with conductor $f\in\Bbb N
$. Then the Dirichlet $L$-series attached to $\chi$ is defined as
$$ L(s,\chi) = \sum_{n=1}^\infty \frac{\chi(n)}{n^s}, \text{ for $ s \in
\Bbb{C},$ cf. [7, 8]}.$$ When $\chi=1$, this is the Riemann zeta
function. In [1], $q$-analogue of $\zeta$-function was defined as
follows: For $h \in \Bbb{Z}, s \in \Bbb{C},$
$$ \zeta_q^{(h)}(s,x) = \frac{1-s+h}{1-s} (q - 1) \sum_{n=0}^\infty
\frac{q^{(n+x)h}}{[n+x]^{s-1}} +  \sum_{n=0}^\infty
\frac{q^{(n+x)h}}{[n+x]^{s}}. \tag 2 $$ Note that
$\zeta_q^{(h)}(s,x)$ is an analytic continuation on $\Bbb{C}$
except for $s=1$ with
$$ \zeta_q^{(h)}(1-m,x) = -\frac{\beta_m^{(h,1)}(x,q)}{m}, \text{ for $ m
\in\Bbb N.$} \tag 3 $$ In [1], we easily see that
$$ \beta_m^{(h,1)}(x,q) = -m \sum_{n=0}^\infty
q^{(n+x)h} [n+x]^{m-1} - (q-1)(m+h)  \sum_{n=0}^\infty q^{(n+x)h}
[n+x]^{m}. \tag 4 $$ It follows from (2) that
$$ \lim_{q \rightarrow 1} \zeta_q^{(h)}(s,x) = \zeta(s,x) =
\sum_{n=0}^\infty \frac{1}{(n+x)^s}. $$ By the meaning of the
$q$-analogue of Dirichlet $L$-series, we consider the following
$L$-series:
$$ L_q^{(h)}(s,\chi) = \frac{1-s+h}{1-s} (q - 1) \sum_{n=1}^\infty
\frac{q^{nh} \chi(n)}{[n]^{s-1}} +  \sum_{n=1}^\infty \frac{q^{nh}
\chi(n)}{[n]^{s}}, \text{ cf. [1], } \tag 5 $$ for $h \in \Bbb{Z},
s \in \Bbb{C}$. It is easy to see that $L_q^{(h)}(s,\chi)$ is an
analytic continuation on $\Bbb{C}$ except for $s=1$. For $m \geq
0,$ the generalized extended $q$-Bernoulli numbers with $\chi$ are
defined as
$$ \beta_{m,\chi}^{(h,k)}(q) = [f]^{m-k} \sum_{i_1, \cdots, i_k = 0}^{f-1}
q^{\overset k \to {\underset {l=1} \to \sum} (h-l+1)i_l}
\beta_{m}^{(h,k)} (\frac{\overset k \to {\underset {l=1} \to \sum}
i_l}{f}, q^f)\left(\overset k \to {\underset {j=1} \to \Pi}
\chi(i_j)\right), \text{ see [1] }. \tag 6 $$  By (5) and (6), we
easily see that
$$ L_q^{(h)}(1-m,\chi)= - \frac{\beta_{m,\chi}^{(h,1)}}{m},\text{ for $m\in \Bbb N$, cf. [1]}. $$

In the present paper we give  twisted $q$-Bernoulli numbers by
using $p$-adic invariant integrals on $\Bbb{Z}_p$. Moreover, we
construct the analogs of $q$-zeta function and $q$-$L$-series
which interpolate the twisted $q$-Bernoulli numbers at negative
integers.


\vskip 0.5 true cm {\bf 2. $q$-extension of Bernoulli numbers}
\par
\vskip 0.5 true cm

In this section we assume that $q \in \Bbb{C}_p$ with $|1- q|_p <
1$. By the definition  of $p$-adic invariant integrals, we see
that
$$ I_1(f_1) = I_1(f) + f^{\prime}(x), \tag 7$$
where $f_1(x) = f(x+1).$ Let
$$T_p = \cup_{n \ge 1}C_{p^n} = \lim_{n\rightarrow \infty}
C_{p^n}, $$ where $C_{p^n} = \{ w | w^{p^n} = 1 \}$ is the cyclic
group of order $p^n$. For $w \in T_p ,$ we denote by $\phi_w :
\Bbb{Z}_p \longrightarrow \Bbb{C}_p$ the locally constant function
$x \mapsto w^x$. If we take $f(x) = \phi_w(x) e^{tx}$, then we
easily see that
$$ \int_{\Bbb{Z}_p} e^{tx} \phi_w(x) d\mu_1(x) = \frac{t}{we^t - 1}, \quad
\text{cf. [5]}. \tag 8 $$ It is obvious from (7) that
$$ \int_{X} e^{tx} \chi(x) \phi_w(x) d\mu_1(x) =
\frac{\overset f \to {\underset {i=1} \to \sum} \chi(i) \phi_w(i)
e^{it}} {w^f e^{ft} - 1}, \text{ cf. [5] }. \tag 9 $$ Now we
define the analogue of Bernoulli numbers as follows:
$$ \aligned
e^{xt}  \frac{t}{we^t - 1} & = \sum_{n=0}^\infty B_{n,w}(x) \frac{t^n}{n!} ,\\
\frac{\overset f \to {\underset {i=1} \to \sum} \chi(i) \phi_w(i)
e^{it}} {w^f e^{ft} - 1} & = \sum_{n=0}^\infty B_{n,w,\chi}
\frac{t^n}{n!}, \text{ cf. [5] }.
\endaligned \tag 10 $$
By (8), (9) and (10), it is not difficult to see that
$$ \int_{\Bbb{Z}_p} x^n \phi_w(x) d\mu_1(x) = B_{n,w} \quad \text{and} \quad
 \int_{X} \chi(x) x^n \phi_w(x) d\mu_1(x) = B_{n,w,\chi} , \tag 11$$
 where $B_{n,w} = B_{n,w}(0)$.

 From (11) we consider twisted $q$-Bernoulli numbers
  using $p$-adic $q$-integral on $\Bbb{Z}_p$. For $w \in T_p$ and $h \in \Bbb{Z},$
 we define the twisted $q$-Bernoulli polynomials as
 $$ \beta^{(h)}_{m,w}(x,q) = \int_{\Bbb{Z}_p} q^{(h-1)y} w^y [x+y]^m d\mu_q(y).
 \tag 12 $$
 Observe that
 $$ \lim_{q\rightarrow 1}  \beta^{(h)}_{m,w}(x, q) = B_{m,w}(x). $$
When $x=0$, we write  $\beta^{(h)}_{m,w}(0,q) =
\beta^{(h)}_{m,w}(q)$, which are called twisted $q$-Bernoulli
numbers. It follows from (12) that
$$ \beta^{(h)}_{m,w}(x,q) = \frac{1}{(1-q)^{m-1}} \sum_{k=0}^m {m \choose k} q^{xk}
(-1)^k \frac{k+h}{1 - q^{h+k} w}. \tag 13 $$ The Eq.(13) is
equivalent to
$$\beta^{(h)}_{m,w}(q) = -m \sum_{n=0}^\infty [n]^{m-1}q^{hn} w^n -
(q-1)(m+h) \sum_{n=0}^\infty [n]^{m}q^{hn} w^n. \tag 14 $$ From
(13), we obtain the below distribution relation for the twisted
$q$-Bernoulli polynomials as follows: For $n \ge 0,$
$$ \beta^{(h)}_{n,w}(x,q) = [f]^{n-1} \sum_{a=0}^{f-1} w^a q^{ha} \beta^{(h)}_{n,w^f}(\frac
{a+x}f,q^f). $$ Let $\chi$ be the Dirichlet character with
conductor $f \in \Bbb N$. Then we define the generalized twisted
$q$-Bernoulli numbers as follows: For $n \ge 0$,
$$ \beta^{(h)}_{m,w,\chi}(q) = \int_X \chi(x) q^{(h-1)x}w^x [x]^m
d\mu_q(x). \tag 15 $$  By (15), it is easy to see that
$$ \beta^{(h)}_{m,w,\chi}(q) = [f]^{m-1} \sum_{a=0}^{f-1} \chi(a) w^a q^{ha}
\beta^{(h)}_{m,w^f}(\frac af,q^f). \tag 16 $$  Remark. We note
that $\lim_{q\rightarrow 1} \beta^{(h)}_{m,w,\chi}(q) =
B_{m,w,\chi},$ ( see Eq. (10) ).


\vskip 0.5 true cm {\bf 3. $q$-zeta functions}
\par
\vskip 0.5 true cm

In this section we assume that $q \in \Bbb{C}$ with $|q| < 1$.
Here we construct the twisted $q$-zeta function
 and the twisted $q$-$L$-series (see Eq.(2) and Eq.(3)). Let $\Bbb{R}$ be the field of
real numbers and let $w$ be the $p^r$-th root of unity. For $q \in
\Bbb{R}$ with $0 < q < 1, s \in \Bbb{C}$ and $h \in \Bbb{Z},$ we
define the twisted $q$-zeta function as
$$ \zeta_{q,w}^{(h)}(s) = \frac{1 - s + h}{1 - s} (q-1)
\sum_{n=1}^\infty \frac{q^{nh} w^n}{[n]^{s-1}} + \sum_{n=1}^\infty
\frac{q^{nh} w^n}{[n]^{s}}. \tag 17 $$ Note that
$\zeta_{q,w}^{(h)}(s)$ is an analytic continuation on $\Bbb{C}$
except for $s=1$ and $\underset {q\rightarrow 1} \to \lim
\zeta^{(h)}_{q,w}(s) =
\zeta(s,w)=\sum_{n=1}^{\infty}\frac{w^n}{n^s},$ cf. [4].
 We see,
by (17), that
$$ \zeta^{(s-1)}_{q,w}(s) = \sum_{n=1}^\infty \frac{q^{n(s-1)}
w^n}{[n]^{s}}. \tag 18 $$

In what follows, the notation $ \zeta^{(s-1)}_{q,w}(s)$ will be
replaced by  $\zeta_{q,w}(s)$, that is,
$$ \zeta_{q,w}(s) (=  \zeta^{(s-1)}_{q,w}(s))= \sum_{n=1}^\infty \frac{q^{n(s-1)}
w^n}{[n]^{s}}. $$ We note that Eq.(18) is the $q$-extension of
Riemann zeta function. By (14) and (17) we give the values of
$\zeta^{(h)}_{q,w}(s)$ at negative integers as follows: For $m \in
\Bbb N,$ we have
$$ \zeta^{(h)}_{q,w}(1-m) = - \frac{\beta^{(h)}_{m,w}(q)}{m}. \tag
19 $$ By (17), we also see that
$$ \zeta_{q,w}(1-m) = \sum_{n=1}^\infty [n]^{m-1} q^{-mn} w^n . \tag 20 $$
The Eq.(20) seems to be the $q$-analogue of  Euler divergence
theorem for Riemann zeta function. Now we also consider the
twisted $q$-analogue of Hurwitz zeta function as follows: For $s
\in \Bbb{C},$ define
$$ \zeta^{(h)}_{q,w}(s,x) = \frac{1 - s + h}{1 - s} (q-1)
\sum_{n=0}^\infty \frac{q^{(n+x)h} w^n}{[n+x]^{s-1}} +
\sum_{n=0}^\infty \frac{q^{(n+x)h} w^n}{[n+x]^{s}}. \tag 21 $$

 Note that $\zeta_{q}^{(h)}(s,x)$ has an analytic
continuation on $\Bbb{C}$ with only one simple poles at  $s=1$. By
Eq.(13), Eq.(14) and Eq.(21), we obtain the following:
$$  \zeta^{(h)}_{q,w}(1-m,x) = - \frac{\beta^{(h)}_{m,w}(x,q)}{m},
\quad \text{for} \ m > 0. $$

Now we consider the twisted $q$-$L$-series which interpolate
twisted generalized $q$-Bernoulli numbers as follows: For $s \in
\Bbb{C},$ define
$$ L^{(h)}_{q,w}(s, \chi) = \frac{1 - s + h}{1 - s} (q-1)
\sum_{n=1}^\infty \frac{q^{nh} w^n \chi(n)}{[n]^{s-1}} +
\sum_{n=1}^\infty \frac{q^{nh} w^n \chi(n)}{[n]^{s}}, \tag 22 $$
where $w$ is the $p^r$-th root of unity.

For any positive integer $m$ we have
$$ L^{(h)}_{q,w}(1-m, \chi) = - \frac{\beta^{(h)}_{m,w,\chi}(q)}{m}.
$$
The Eq.(22) implies that
$$ \aligned
L^{(s-1)}_{q,w}(s, \chi) & = \sum_{n=1}^\infty \frac{q^{n(s-1)}
w^n
\chi(n)}{[n]^{s}} \\
& = [f]^{-s} \sum_{a=1}^{f} \chi(a) w^a q^{(s-1)a}
\zeta_{q^f,w^f}(s,\frac af). \endaligned $$

\noindent {\bf Question.} Find a $q$-analogue of the $p$-adic
twisted $L$-function which interpolates $q$-Bernoulli numbers
$\beta^{(h)}_{m,w,\chi}(q)$, cf. [3].

ACKNOWLEDGEMENTS: This paper was supported by Korea Research
Foundation Grant( KRF-2002-050-C00001).


\vskip0.5cm

\noindent Institute of Science Education, Kongju National
University, Kongju 314-701, Korea, tkim\@kongju.ac.kr \par
\noindent Department of Mathematics and Computer Science, KonKuk
University, Choongju, Chungbuk 380-701, Korea,
leechae.jang\@kku.ac.kr \par \noindent Department of Mathematics
Education, Kyungpook University, Taegu,
 702-701,  Korea, shrim\@kyungpook.ac.kr\par
 \noindent Faculty of Informatoin and
Science, Daegu Haany University, Kyungsan,  712-240, Korea,
hkpak\@ik.ac.kr

\enddocument